# Euler's Graph World
## - Purity, Regularity and Evenness
## -Law of Nature?
## - Constructions and Examples

Suryaprakash Nagoji Rao*
Surusha5152@hotmail.com
suryaprakashnrao@iitbombay.org<sup>∞</sup>

**ABSTRACT.** We propose a Law of Nature? *Viz., 'Pure Regularity Occurs at Naïve Levels and Regularity has Affinity with Evenness'*. In a series of three papers, it was established that regular Euler graphs with only one type of (pure) cycles are nonexistent; Regular Euler graphs with only two types of cycles are possible in one of the six cases, viz., regular bipartite Euler graphs of degree >2; Evenness plays role in unveiling regularity; Lastly, $K_5$ is a regular Euler graph with three types of cycles (0,1,3); This is the only known graph with the property; It is conjectured that regular Euler graphs of order >5 with only three cycle types are nonexistent and this is proved true in part cases in each of the four cases. Some constructions and examples are given for the Euler graphs under (mod 4) satisfying intersection (combined cycle) rules. These infinite classes of Euler graphs serve as candidates for gracefulness. Infinite families of graceful graphs are presented in Case-0.

## INTRODUCTION

Graph Theory, a branch of Discrete Mathematics, had its birth in different circumstances; a popular context being Konigsberg Bridge Problem of 17th Century. Leonhard Euler, a pioneering Swiss mathematician and physicist, not only answered the non-existence of a closed trail but also Euler's work lead to what are today known as Euler graphs, see Euler (1736,1758), Joli, et al (2010). We refer to Harary (1972) for terminology and notation. Here, by graph we mean a finite, undirected, connected graph without loops and multiple edges.

A *(p,q)-graph* has p nodes and q edges; p and q are *order* and *size* of the graph respectively. A *cycle graph* is a graph that consists of a single cycle. A cycle graph of length n>2 is denoted by $C_n$ or *n-cycle*. A cycle graph is even or odd according as n is even or odd respectively. *Euler trail* is a trail in a graph which visits every edge exactly once. Similarly, Euler circuit is Euler trail which starts and ends on the same node. *Euler graph* admits Euler circuits and are characterized by the simple criterion that all nodes are of even degree. That is, evenness of node degrees characterizes Euler graphs. Euler graphs exhibit many interesting properties. E.g., Edge set of Euler graph may be decomposed into edge disjoint cycles and is called cycle decomposition (CD). Euler graph may have two or more CDs. Euler graph and a CD is shown in Fig.1a and Fig.1b,c. Here, Euler graph is decomposed into five edge disjoint cycles: Three 4-cycles and two 6-cycles.

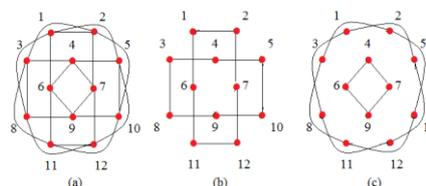

Fig.1 Euler graph (a) and a cycle decomposition (b&c).

A regular graph has all node degrees same. The simplest regular graphs are the cycle graphs and are regular Euler graphs (See Fig.2a,b,c,d for the first four). On the other extreme the complete graph $K_p$ is regular of degree p-1 and for p>5 has all four cycle types $C_n$, n≡0,1,2or3(mod 4). $K_p$ is Euler for every p≥3 odd. Triangle is the first regular Euler graph with the only cycle of length 3. Cycle graph $C_4$ is the only regular Euler graph of order 4. Complete graph of order 4, $K_4$ is the other regular graph of degree three having only cycle types $C_n$, n≡0&3(mod 4) and not Euler. For order 5 there are two regular Euler graphs, viz., $C_5$ and $K_5$ and for order 6, $C_6$ and the graph in Fig.2g are regular Euler of degree 2 and 4 respectively. Fig.2f has only 0,1,3 cycle types. Figs.2g,h have all four cycle types.

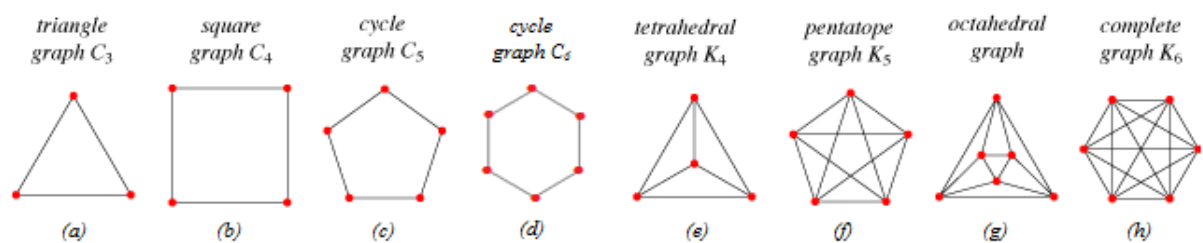

Fig.2 Regular Graphs.

This implies that regularity within structurally simple Euler graph systems is restricted and realized at low levels. By level here we mean graph size or degree of regularity of the graph. Regular bipartite graphs of degree larger than two are sparse with reference to order. For example, Hypercube or n-Cube is a regular bipartite (p,q)-graph of degree n for n≥0 and exists for $p=2^n$ only. Hypercube is regular Euler for n even. A regular bipartite Euler graph viz., 4-cube is shown in Fig.3.
http://en.wikipedia.org/wiki/Hypercube.

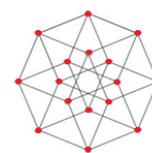

Fig.3 4-Cube

A graph G is called a *'labeled graph'* when each node u is assigned a label φ(u) and each edge uv is assigned the label φ(uv)=|φ(u)-φ(v)|. In this case φ is called a *'labeling'* of G. Define $N(φ)=\{n∈\{0,1,...,q_0\}: φ(u)=n,$ for some $u∈V\}$, $E(φ)=\{e∈\{1,2,...,q_0\}: |φ(u)-φ(v)|=e,$ for some edge $uv∈E(G)\}$. Elements of N(φ), (E(φ)) are called *'node (edge) labels'* of G with respect to φ. A (p,q)-graph G is *'gracefully labeled'* if there is a labeling φ of G such that N(φ)⊆{0,1,...,q} and E(φ)={1,2,...,q}. Such a labeling is called a *'graceful labeling'* of G. A *'graceful graph'* can be gracefully labeled, otherwise it is a 'nongraceful graph' (see Rosa (1967), Golomb (1972). Varied applications of labeled graphs have been cited in the literature survey by Gallion (2013).

**EULER'S GRAPH WORLD** (Jun-July 2014)

We consider the system of *Euler's Graph World* under the (mod 4) operation. This divides the class of Euler graphs into four classes viz., Euler graphs with only cycle types $C_n$, n≡0,1,2,3(mod 4). Denote the class of graphs with exactly j cycle types by $\mathcal{C}_j$, j=1,2,3,4. That is, $\mathcal{C}_j$ is the class of Euler graphs with j number of cycle types out of the four ≡0,1,2,3(mod 4). For example, $\mathcal{C}_1$ is the class of Euler graphs with only one type of cycles $C_n$, n≡i(mod 4), i=0,1,2,3. Further, denote by $\varepsilon_i$, $\varepsilon_{ij}$, $\varepsilon_{ijk}$, (i,j,k=0,1,2,3) the class of Euler graphs with exactly one, two, three types of cycles. $\varepsilon_{0123}$ is the class of Euler graphs with all four types of cycles.

**Class $\mathcal{C}_1$**

This is the class of Euler graphs with only one type of cycles (See Rao (2014), Part-I and Rao (2020)). We say in this case Euler graph has *pure* cycles, which means every cycle decomposition of Euler graph in $\mathcal{C}_1$ has exactly one class in any of its cycle decomposition. Otherwise, the Euler graph is said to have impurity or has cycles of at least two of the four types under (mod 4). It is established that

**Theorem** 1. *The class of Euler graphs $\mathcal{C}_1$, with only one type of cycles under (mod 4) operation, is regular only when they are cycle graphs.*

So, cycle graphs, regular of degree 2, are the only regular graphs with only one type of cycles $C_n$, n≡i(mod 4) for i=0,1,2,3 allowed. This characterizes regular Euler graphs of class $\mathcal{C}_1$:

**Theorem** 2. *Euler graphs with only one type of cycles are regular iff they are cycle graphs.*

Euler graphs need not be bipartite, so a corollary follows that

**Theorem** 3. *Euler graphs with only one type of cycles are regular bipartite iff they are even cycle graphs.*

Thus in summary, Euler graphs under (mod 4) operation exhibit the property viz.,

**Property** 1. *Regularity amongst Euler graphs with only one type of cycles is realized at naïve levels.*

Thus we infer that

**Property** 2. *Purity manifests only at naive levels of regularity.*

Further,

**Theorem** 4. *Euler graphs in $\varepsilon_i$, i=1,2,3 are planar. Euler graphs in $\varepsilon_0$ need not be planar.*

Examples of nonplanar graphs in $\varepsilon_0$ are given.

More information on regularity and evenness is given in Appendix.

**Class $\mathcal{C}_2$**
This is the class of Euler graphs with only two types of cycles $C_n$ or cycle decomposition contains exactly two classes (See Rao (2014), Part-II and Rao (2020)). This gives raise to six cases, viz., n≡i&j(mod 4) for (i,j)=(0,1), (0,2), (0,3), (1,2), (1,3), (2,3).

It is established that

*Theorem 5. The class of Euler graphs $\mathcal{C}_2$ is devoid of regularity except when all the cycles $C_n$ of Euler graph are of even length, that is, when (i,j)= (0,2) or n≡0&2(mod 4).*

These are precisely bipartite Euler graphs having exactly two classes in any of its cyclic decompositions. In other words, regular bipartite Euler graphs of degree >2 are the only Euler graphs from the class $\mathcal{C}_2$ and exhibit regular structures of different complexity for higher orders and higher degree of regularity. In the other five cases, viz., n≡i&j(mod 4) for (i,j)=(0,1), (0,3), (1,2), (1,3), (2,3) regularity is non-existent. This characterizes $\mathcal{C}_2$ as follows:

*Theorem 6. Euler graph with only two types of cycles is regular iff it is regular bipartite Euler graph of degree >2.*

We summarize the results for the Classes $\mathcal{C}_1$ and $\mathcal{C}_2$ below:

*Theorem 7. Euler graph with utmost two types of cycles is regular iff it is either an odd cycle graph or a regular bipartite Euler graph.*

Further,

*Theorem 8. Euler graphs in $\varepsilon_{ij}$, (i,j)=(1,2),(1,3),(2,3) are planar.*

Examples of nonplanar graphs in $\varepsilon_{01}, \varepsilon_{02}, \varepsilon_{03}$ are given.

We infer that

*Property 3. Complex regularity in Euler graphs is manifested in presence of impurity or mixture of at least two types of cycles.*

This indicates that

*Property 4. There is a minimum level of impurity for complex regularity to exhibit or unveil or to show up.*

Evenness also has its significant role to play in manifestation of the regular structures. Higher degree regularity unveils in the graphs of class $\mathcal{C}_2$ when not only the node degrees are even but also all the cycles in a graph are of even length.

**Class $\mathcal{C}_3$**
This is the class of Euler graphs with only three types of cycles or graphs having exactly three of the four classes in any cycle decomposition (See Rao (2014), Part-III). There are four cases with only three cycle types, $C_n$, viz., n≡i,j&k(mod 4) for (i,j,k)=(0,1,2), (0,1,3), (0,2,3), (1,2,3) under the (mod 4) operation. Such a graph with three distinct types of cycles i,j,k from {0,1,2,3} belongs to $\varepsilon_{ijk}$.

An example of a regular graph of degree 4 with only cycle types (0,1,3) is $K_5$, see Fig.2f. That cycles of the type $C_n$, n≡2(mod 4) are not in $K_5$ is trivial. $K_5$ is not only the smallest graph with this property but also $K_5$ is the only regular Euler graph so far known from $\mathcal{C}_3$ in all the above four cases. With no success in several attempts to prove or search or construct such graphs we propose:

**Conjecture**. *Regular Euler graphs of order p>5 with exactly three cycle types are nonexistent.*
*or*
*Regular Euler graphs of order p>5 with three cycle types also contains fourth type.*

Even if the Conjecture is wrong, regular Euler graphs in $\mathcal{C}_3$ is supposedly sparsely distributed.

*Theorem 5. Euler graphs in $\varepsilon_{ijk}$, (i,j,k)=(0,1,2), (0,1,3), (0,2,3) may be nonplanar.*

A nonplanar graph is given in each case. Nonplanarity in the case (1,2,3) is not known.

**Class $\mathcal{C}_4$**

This is the class of Euler graphs with all four types of cycles present or every cycle decomposition has exactly four classes. Most of the regular Euler graphs are of this type and the graphs are not necessarily planar.

The above inferences in the classes $\mathcal{C}_j$, j=1,2,3,4 guide to propose:

**Law of Nature?**

- In Systems of Nature
    a. *Purity and Regularity don't go together except at naïve levels.*
    b. *Regularity has affinity with Evenness.*

**Constructions and Examples** (Jun-Aug 2020)

We give infinite families of graphs in each case. These serve as candidates for infinite families of graphs for gracefulness.

**Euler Graphs with Only One Type of Cycles**

| Table-1 Case 0,1,2,3 | | | |
|---|---|---|---|
| Cycle Type | | Combined Cycle type | |
| Cycle-1 | Cycle-2 | Int Even | Int Odd |
| 0 | 0 | 0 | 2 |
| 1 | 1 | 2 | 0 |
| 2 | 2 | 0 | 2 |
| 3 | 3 | 2 | 0 |

Cycles in any cycle decomposition of Euler graph are only one type. In this section we consider Euler graphs with only one type of cycles are considered under the operation (mod 4) (See Rao (2014), [4]). Cycle intersection (Combined cycle) rules for Euler graphs from $\varepsilon_i$, i=0,1,2,3 are given in the Table-1.

**Case-0. $\varepsilon_0$: Euler graphs with only cycles $C_n$, n≡0(mod 4)**

*'Cycle intersection (Combined cycle) rules'* for Euler graphs from $\varepsilon_0$ are given in the Table-1. Two cycles under ≡0(mod 4) intersect in even path. Two cycles under ≡0(mod 4) intersect in odd path lead to cycle of type ≡2(mod 4) which is shown in yellow. Similarly, two cycles under ≡2(mod 4) intersect in odd path. Two cycles under ≡1 or 3(mod 4) cannot intersect and hence Euler graphs from $\varepsilon_i$, i=1,3 are either cycle graphs or each block of the graph is a cycle.

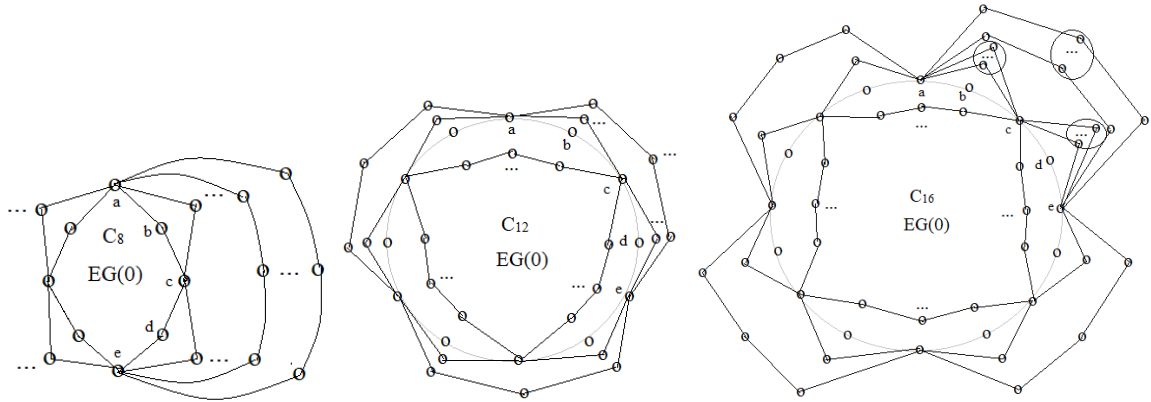

Fig. 4 Infinite classes of graphs from $\varepsilon_0$.

Handle at a u-v path, H(u,v) is defined using an example. For a cycle $C_{16}$ and the a-e path a,b,c,d,e (see Fig.4c) a *'handle'* H(a,e) at a,e is defined as a collection of paths originating at a and ending at e of length 4 containing only a,e nodes of the cycle (such paths have even intersection follows from Table-1). Number of paths added is the *size* of the handle. In this case, paths of length 2 handles H(a,c) and H(c,e) also may be considered at a,c and c,e.

Note that, adding a path of length m at a u-v path of given cycle with d(u,v)=m has the property that cycles containing this new u-v path creates cycle of the type already existing in G. Hence the resulting graph with handle added is of the same type as G. The handle may contain any number of u-v paths making the class of graphs infinite. If the starting cycle is of large length handles of different type may be added. In this case, handles of length 2 and 4 are considered (see Fig.4c). Same explanation holds for Fig.4a,b. Three dots indicate that the handle may be considered with repeated path. Same convention is followed throughout the paper.

**Infinite Family of Graceful Euler Graphs From $\varepsilon_0$:**

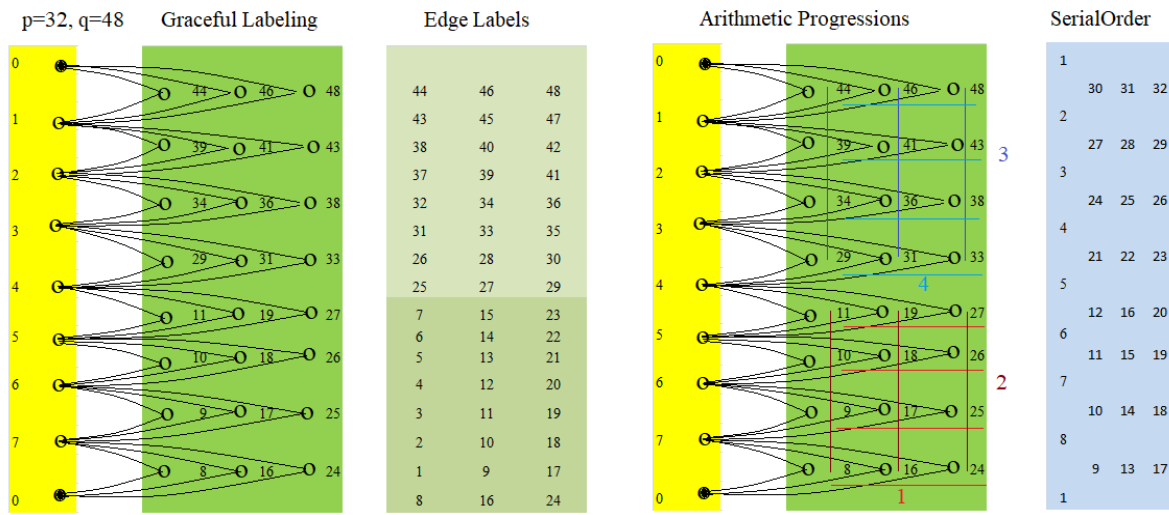

Fig. 4(i) Euler graph with graceful labeling and edge labels. Arithmetic progressions.

The (32,48)-graph in Fig.4(i)a is same as the (32,48)-graph in Fig.4(ii). It is an example of Euler graph with only one type of cycles, in this case cycles $\equiv 0 \pmod 4$. It has a 16-cycle with a K(2,3) on the alternate nodes as shown. A graceful labeling is also shown. Edge labels arising out of the graceful labeling are shown in Fig.4(i)c. The graceful labeling of the graph in Fig.4(i)c consists of arithmetic progressions. Red color horizontal lines in the graph Fig.4(i)c indicate the AP1 (8,16,24) and has first term a=2t, common difference d=2t. Brown color vertical lines indicate the AP2 (8,9,10,11) and has a=2t+1, d=1. Navy blue color vertical lines indicate the AP3 (29,34,39,44) and has a=2st+(t-1)+2, d=2s-1. Blue color horizontal lines indicate AP4 (29,31,33) and has a=2st+(t-1)+2, d=2. The serial order of graceful numbering followed is shown in Fig.4(i)d.

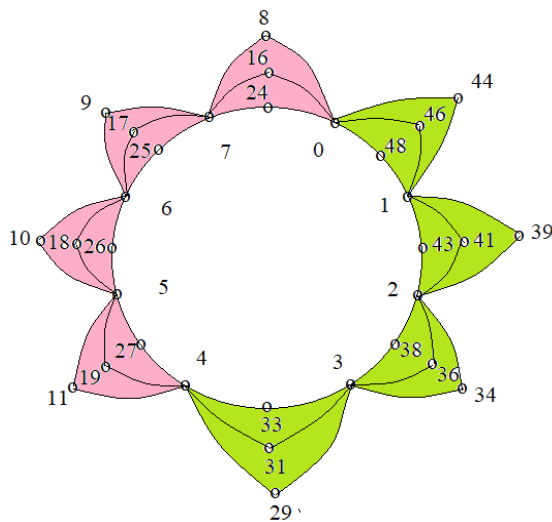

We define graph G(t,s) consisting of a cycle $C_{4t}$, t≥1 and a K(2,s) s≥1 on the alternate nodes of the cycle $C_{4t}$ as shown in Fig4(i) and (ii) graph G(4,3) for t=4 and s=3. Order and size of G(t,s) are p=2t(s+1) and q=4ts. Graceful numbering of this graph is given in Table-1a below. Serial order of the nodes for the graceful numbering is given in Table-1b. The graceful numbering consists of arithmetic progressions horizontally and vertically as illustrated above for the graph in Fig.4(i).

Fig. 4(ii) Same graph as in Fig.4(i) with graceful labeling.

Table-1a. Graceful Numbering

| | | | | | | |
|---|---|---|---|---|---|---|
| 0 | 2st+(t-1)1+2+(t-1)(2s-1) | 2st+(t-1)1+2+(t-1)(2s-1)+2 | 2st+(t-1)1+2+(t-1)(2s-1)+(l-1)2 | | 2st+(t-1)1+2+(t-1)(2s-1)+(s-1)2 | |
| … | … | … | … | | … | |
| | 2st+(t-1)1+2+(k-1)(2s-1) | 2st+(t-1)1+2+(k-1)(2s-1)+2 | 2st+(t-1)1+2+(k-1)(2s-1)+(l-1)2 | | 2st+(t-1)1+2+(k-1)(2s-1)+(s-1)2 | |
| | … | … | … | | … | |
| t-2 | 2st+(t-1)1+2+(2s-1) | 2st+(t-1)1+2+(2s-1)+2 | 2st+(t-1)1+2+(2s-1)+(l-1)2 | | 2st+(t-1)1+2+(2s-1)+(s-1)2 | |
| t-1 | 2st+(t-1)1+2 | 2st+(t-1)1+2+2 | … | 2st+(t-1)1+2+(s-1)2 | … | 2st+(t-1)1+2+(s-1)2 |
| t | 2t+(t-1)1=3t-1 | 4t+(t-1)1=5t-1 | 2it+(t-1)1=(2i+1)t-1 | | 2st+(t-1)=(2s+1)t-1 | |
| … | … | … | … | | … | |
| 2t-3 | 2t+2 | 4t+2 | 2it+2 | | 2st+2 | |
| 2t-2 | 2t+1 | 4t+1 | 2it+1 | | 2st+1 | |
| 2t-1 | 2t | 4t | … | 2it | … | 2st |

Table-1b. Serial order of Graceful Numbering

| | | | | | |
|---|---|---|---|---|---|
| 1 | t(s+2)+1+(t-1)+(s-1)(t-1) | t(s+2)+1+(t-1)+(s-1)(t-1)+1 | … | | t(s+2)+1+(t-1)+(s-1)(t-1)+(s-1)1 |
| … | … | … | … | | … |
| t-1 | t(s+2)+1+(s-1)2+1 | t(s+2)+1+(s-1)2+1+1 | | | t(s+2)+1+(s-1)2+(s-1)1 |
| t | st+2t+1=t(s+2)+1 | t(s+2)+1+1 | | | t(s+2)+1+(s-1)1 |
| t+1 | 2t+t=3t | 3t+t=4t | | | (s+1)t+1+(t-1)1=st+2t |
| … | … | … | … | | … |
| 2t-1 | 2t+2 | 3t+1+1 | | | (s+1)t+1+1 |
| 2t | 2t+1 | 3t+1 | | | 2t+1+(s-1)t=(s+1)t+1 |

**Truncated Infinite Families of Graceful Euler Graphs From ε₀:**

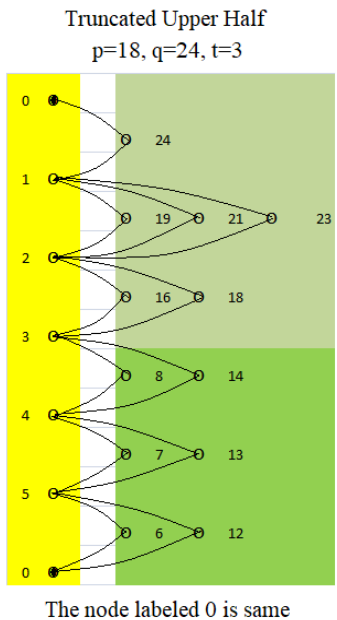

Truncated Upper Half
p=18, q=24, t=3

The node labeled 0 is same

The above family of graceful Euler graphs can be modified by truncating the upper half portion resulting in several infinite families of graceful Euler graphs. The same numbering rule works as a graceful numbering. We illustrate this by an example (see Fig.4(iii)): p=24 , q=36 , t=3. Note that the graph is not Euler. Truncation of the upper half can be done in many ways. Each row in light green portion shall have at least one node. Each row in green portion shall have equal nodes making it aqn array.

Last row in the green portion nodes are labeled 2t, 4t,… The t-1 nodes above them are labeled by the APs with a=2t+1,d=1; … First node in the last row of light green portion is numbered 2 more than the node number of the last term (top) in the last column of green portion. Last row in light green portion is numbered by the AP with d=2. The next row in the light green portion is AP with d=2. First node is one more than the last node number of the previous row. See also the graceful numbering in Fig.4(i).

No details are given here.

Fig. 4(iii) Truncated Graceful Family.

**Case-1. $\varepsilon_1$: Euler graphs with only cycles $C_n$, n≡1(mod 4)**

Cycle intersection (Combined cycle) rules for Euler graphs from $\varepsilon_1$ are given in the Table-1. Graphs in this class are simple and have the structure with each block a cycle graph $C_i$, i≡1(mod 4).

**Case-2. $\varepsilon_2$: Euler graphs with only cycles $C_n$, n≡2(mod 4)**

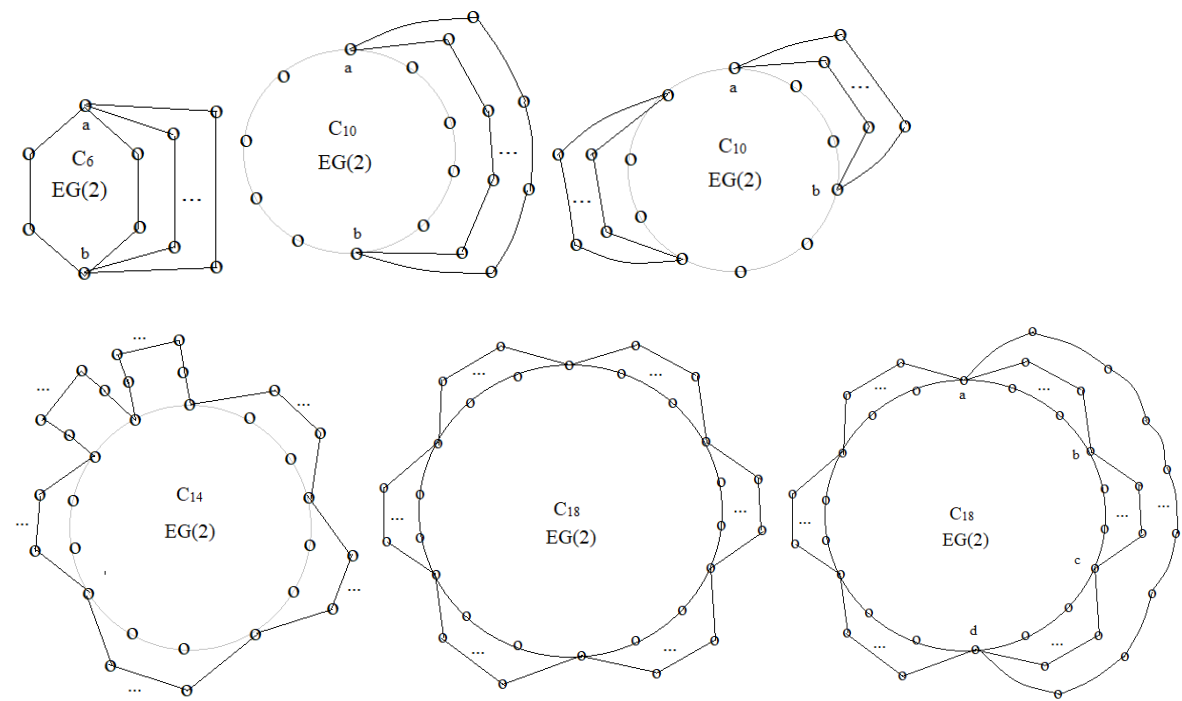

Fig. 5 Infinite classes of graphs from $\varepsilon_2$.

Cycle intersection (Combined cycle) rules for Euler graphs from $\varepsilon_2$ are given in the Table-1. In the case, Fig.5d there are four handles of length 3 and two handles of length 6. The three dots at a handle H(u,v) indicate that any number of u-v paths may be added as the new cycles are of the same type that of G. The handle size (number of paths) is so chosen to make the resulting graph Euler. Fig.5e has six handles and Fig.5f has additionally a repeatable handle of length 9.

**Case-3. $\varepsilon_3$: Euler graphs with only cycles $C_n$, n≡3(mod 4)**

Cycle intersection (Combined cycle) rules for Euler graphs from $\varepsilon_3$ are given in the Table-1. Graphs in this class are simple and have the structure with each block a cycle graph $C_i$, i≡3(mod 4).

**Euler Graphs with Only Two Types of Cycles**

A natural generalization of the class of trees viz., graphs with only one type of cycles under (mod 4) operation was studied in Rao (2014) [4]. Here, we consider the class of graphs with only two types of cycles in any cycle decomposition of Euler graph (See Rao (2014) [5]). Denote the class of Euler graphs G by $\varepsilon_{ij}$ with only two classes $\mathcal{C}_i$ and $\mathcal{C}_j$ in every cycle decomposition of G, that is, graphs having only two types of cycles $C_i$ and $C_j$, i≠j, i,j≡0,1,2 or 3(mod 4). This assumption simplifies considerably the general structure of Euler graphs and stand next to trees and graphs with only one type of cycles including unicyclic graphs.

Six cases arise with only two types of cycles under (mod 4) operation: 1) n≡0&1(mod 4), 2) n≡0&2(mod 4), 3) n≡0&3(mod 4), 4) n≡1&2(mod 4), 5) n≡1&3(mod 4), 6) n≡2&3(mod 4). Case-2 corresponds to the class of bipartite Euler graphs.

**Case-01. $\varepsilon_{01}$: Euler graphs with only cycles $C_n$, n≡0&1(mod 4)**
(See Rao (2014), [5])

| Table-2 | Case 0,1 | | |
|---|---|---|---|
| Cycle Type | | Combined Cycle type | |
| Cycle-1 | Cycle-2 | Int Even | Int Odd |
| 0 | 0 | 0 | 2 |
| 0 | 1 | 1 | 3 |
| 1 | 1 | 2 | 0 |

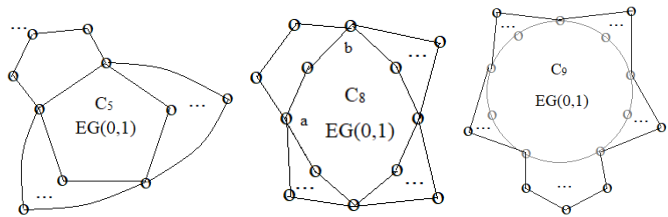

Fig. 6 Infinite class of graphs from $\varepsilon_{01}$.

Cycle intersection (Combined cycle) rules for Euler graphs from $\varepsilon_{01}$ are given in the Table-2. In Fig.6b the Handle(a,b) is a non-repeatable handle.

**Case-02. $\varepsilon_{02}$: Euler graphs with only cycles $C_n$, n≡0&2(mod 4)**

| Table-3 | Case 0,2 | | |
|---|---|---|---|
| Cycle Type | | Combined Cycle type | |
| Cycle-1 | Cycle-2 | Int Even | Int Odd |
| 0 | 0 | 0 | 2 |
| 0 | 2 | 2 | 0 |
| 2 | 2 | 0 | 2 |

Cycle intersection (Combined cycle) rules for Euler graphs from $\varepsilon_{02}$ are given in the Table-3. Note that the intersection of cycles ≡0 or 2 result in type 0 or 2. This is the well known class of bipartite graphs.

**Case-03. $\varepsilon_{03}$: Euler graphs with only cycles $C_n$, n≡0&3(mod 4)**

| Table-4 | Case 0,3 | | |
|---|---|---|---|
| Cycle Type | | Combined Cycle type | |
| Cycle-1 | Cycle-2 | Int Even | Int Odd |
| 0 | 0 | 0 | 2 |
| 0 | 3 | 3 | 1 |
| 3 | 3 | 2 | 0 |

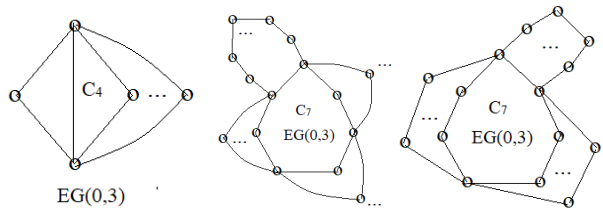

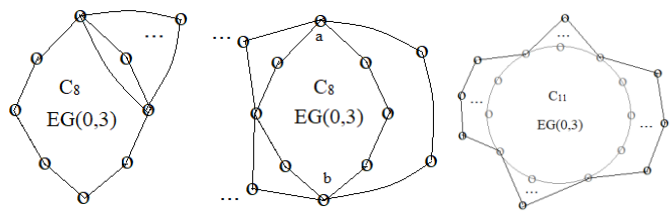

Fig. 7 Infinite classes of graphs from $\varepsilon_{03}$.

Cycle intersection (Combined cycle) rules for Euler graphs from $\varepsilon_{03}$ are given in the Table-4. Fig.7e has a handle H(a,b) of length 3 and cannot be repeated (no dots) else there exists a cycle $\equiv 2(\bmod 4)$.

**Case-12. $\varepsilon_{12}$: Euler graphs with only cycles $C_n$, n≡1&2(mod 4)**

| Table-5 | Case 1,2 | | |
|---|---|---|---|
| Cycle Type | | Combined Cycle type | |
| Cycle-1 | Cycle-2 | Int Even | Int Odd |
| 1 | 1 | 2 | 0 |
| 1 | 2 | 3 | 1 |
| 2 | 2 | 0 | 2 |

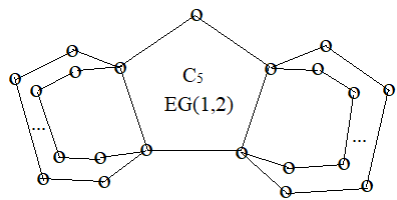

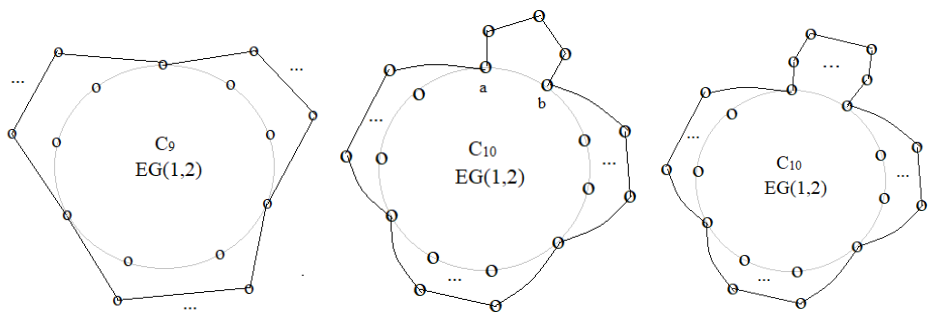

Fig. 8 Infinite classes of graphs from $\varepsilon_{12}$.

Cycle intersection (Combined cycle) rules for Euler graphs from $\varepsilon_{12}$ are given in the Table-5. In Fig.8c, handle H(a,b) of length 4 cannot be repeated else cycle $\equiv 0(\bmod 4)$ results.

**Case-13. $\varepsilon_{13}$: Euler graphs with only cycles $C_n$, n≡1&3(mod 4)**

| Table-6 | Case 1,3 | | |
|---|---|---|---|
| Cycle Type | | Combined Cycle type | |
| Cycle-1 | Cycle-2 | Int Even | Int Odd |
| 1 | 1 | 2 | 0 |
| 1 | 3 | 0 | 2 |
| 3 | 3 | 2 | 0 |

Cycle intersection (Combined cycle) rules for Euler graphs from $\varepsilon_{13}$ are given in the Table-6. This class of graphs consists of each block a cycle graph $\equiv 1$ or $3(\bmod 4)$.

**Case-23. $\varepsilon_{23}$: Euler graphs with only cycles $C_n$, n≡2&3(mod 4)**

| Table-7 | Case 2,3 | | |
|---|---|---|---|
| Cycle Type | | Combined Cycle type | |
| Cycle-1 | Cycle-2 | Int Even | Int Odd |
| 2 | 2 | 0 | 2 |
| 2 | 3 | 1 | 3 |
| 3 | 3 | 2 | 0 |

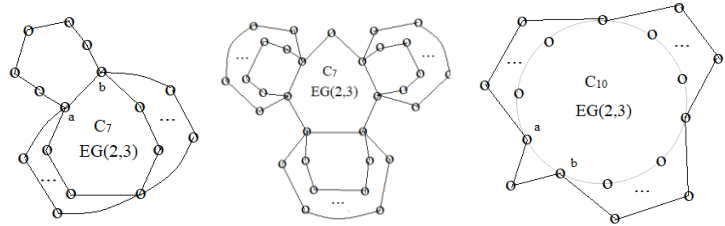

Fig. 9 Infinite classes of graphs from $\varepsilon_{23}$.

Cycle intersection (Combined cycle) rules for Euler graphs from $\varepsilon_{23}$ are given in the Table-7. Graph in Fig.9a has a handle H(a,b) of length 6 and cannot be repeated else a cycle ≡0(mod 4) results. Fig.9b has three handles of same type. The size, number of handles may be reduced to two or one. Fig.9c has three repeatable handles and one handle H(a,b) which cannot be repeated.

**Euler Graphs with Only Three Types of Cycles**

Four cases arise (See Rao (2014), [6]) according as Euler graph has only three of the four cycle types: (0,1,2), (0,1,3), (0,2,3), (1,2,3). A graph with each block being a i-cycle graph, i=0,1,2,3 is Euler. Such a graph with three distinct types of cycles i,j,k from {0,1,2,3} belongs to $\varepsilon_{ijk}$.

**Case-012. $\varepsilon_{012}$: Euler graphs with only cycles $C_n$, n≡0,1&2(mod 4)**

| Table-8 | Case 0,1,2 | | |
|---|---|---|---|
| Cycle Type | | Combined Cycle type | |
| Cycle-1 | Cycle-2 | Int Even | Int Odd |
| 0 | 0 | 0 | 2 |
| 0 | 1 | 1 | 3 |
| 0 | 2 | 2 | 0 |
| 1 | 1 | 2 | 0 |
| 1 | 2 | 3 | 1 |
| 2 | 2 | 0 | 2 |

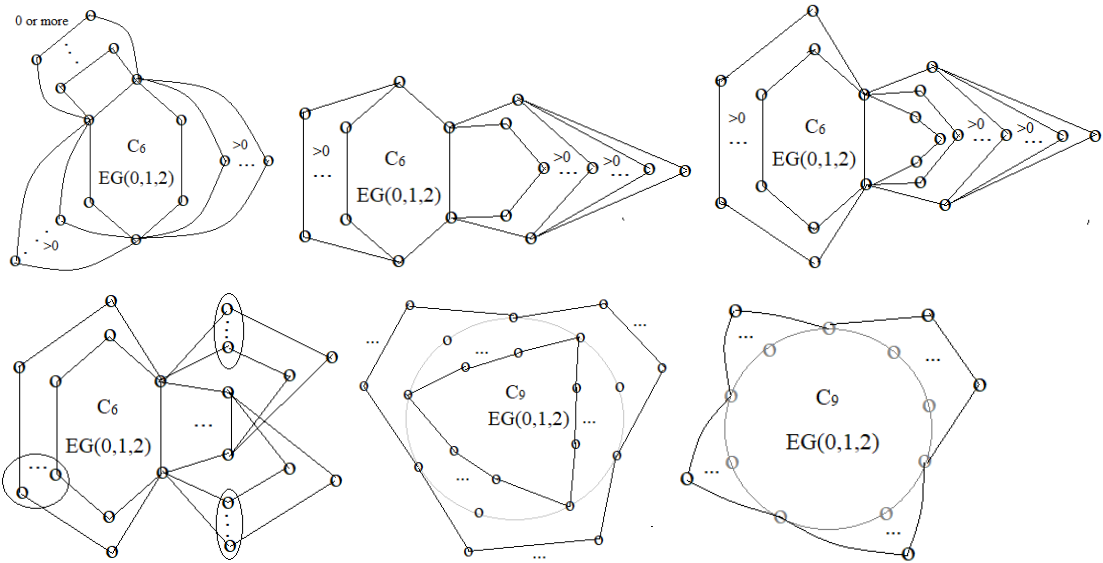

Fig.10 Infinite classes of graphs from $\varepsilon_{012}$.

Cycle intersection (Combined cycle) rules for Euler graphs from $\varepsilon_{012}$ are given in the Table-8. All six graphs have repeatable handles.

Euler's Graph World: Purity, Regularity and Evenness -Law of Nature? Constructions and Examples

**Case-013. $\varepsilon_{013}$: Euler graphs with only cycles $C_n$, n≡0,1&3(mod 4)**

| Table-9 | Case 0,1,3 | | |
|---|---|---|---|
| Cycle Type | | Combined Cycle type | |
| Cycle-1 | Cycle-2 | Int Even | Int Odd |
| 0 | 0 | 0 | 2 |
| 0 | 1 | 1 | 3 |
| 0 | 3 | 3 | 1 |
| 1 | 1 | 2 | 0 |
| 1 | 3 | 0 | 2 |
| 3 | 3 | 2 | 0 |

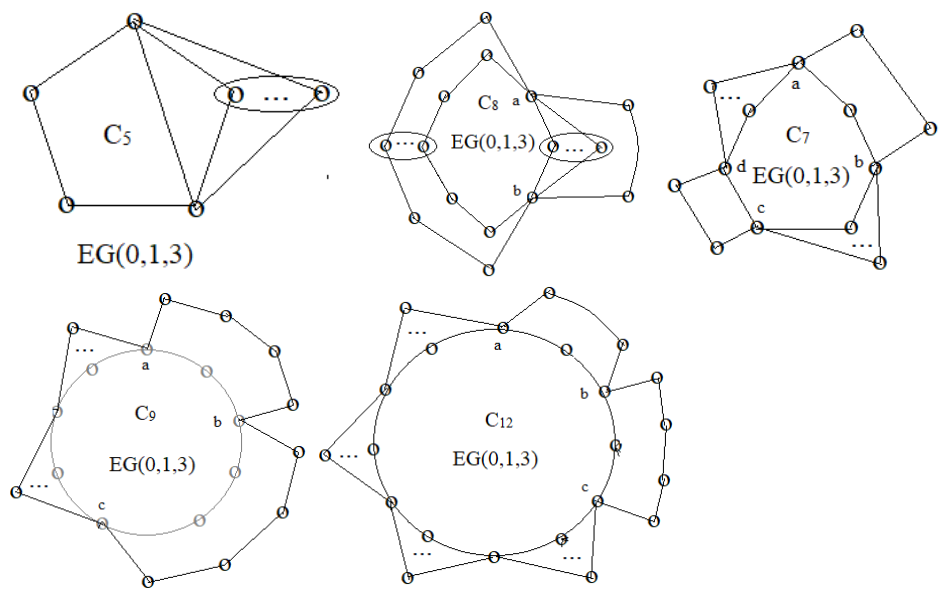

Fig.11 Infinite classes of graphs from $\varepsilon_{013}$.

Cycle intersection (Combined cycle) rules for Euler graphs from $\varepsilon_{013}$ are given in the Table-9. In Fig.11b handles are so chosen to make the graph Euler. Further, handle H(a,b) of length 3 is non-repeatable else it will have a cycle ≡2(mod 4). Fig.11c has non-repeatable handles H(a,b) and H(c,d). Fig.11d has non-repeatable handles H(a,b) and H(b,c) of length 5. Fig.11e has non-repeatable handles H(a,b) and H(b,c) of length 3 and 5 respectively.

**Case-023. $\varepsilon_{023}$: Euler graphs with only cycles $C_n$, n≡0,2&3(mod 4)**

| Table-10 | Case 0,2,3 | | |
|---|---|---|---|
| Cycle Type | | Combined Cycle type | |
| Cycle-1 | Cycle-2 | Int Even | Int Odd |
| 0 | 0 | 0 | 2 |
| 0 | 2 | 2 | 0 |
| 0 | 3 | 3 | 1 |
| 2 | 2 | 0 | 2 |
| 2 | 3 | 1 | 3 |
| 3 | 3 | 2 | 0 |

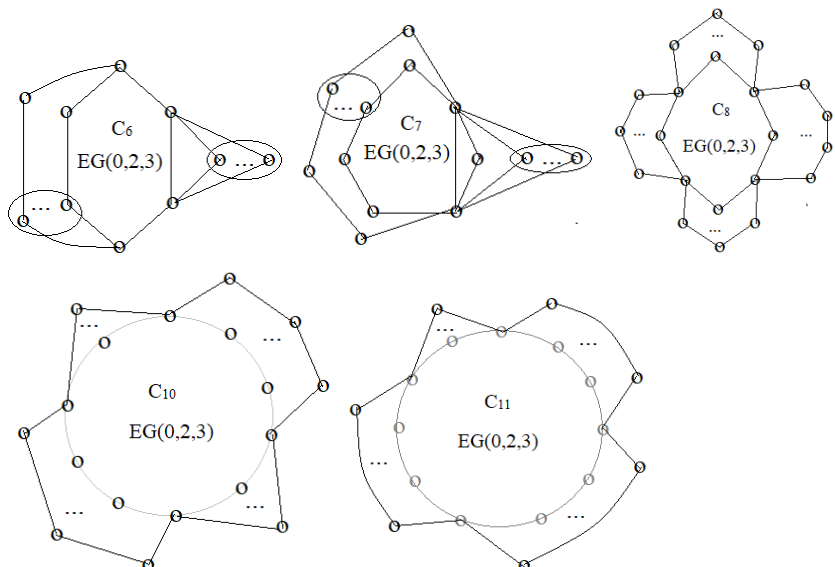

Fig.12 Infinite classes of graphs from $\varepsilon_{023}$.

Cycle intersection (Combined cycle) rules for Euler graphs from $\varepsilon_{023}$ are given in the Table-10. All the five graphs have repeatable handles.

**Case-123. $\varepsilon_{123}$: Euler graphs with only cycles $C_n$, n≡1,2&3(mod 4)**

| Table-11 | Case 1,2,3 | | |
|---|---|---|---|
| Cycle Type | | Combined Cycle type | |
| Cycle-1 | Cycle-2 | Int Even | Int Odd |
| 1 | 1 | 2 | 0 |
| 1 | 2 | 3 | 1 |
| 1 | 3 | 0 | 2 |
| 2 | 2 | 0 | 2 |
| 2 | 3 | 1 | 3 |
| 3 | 3 | 2 | 0 |

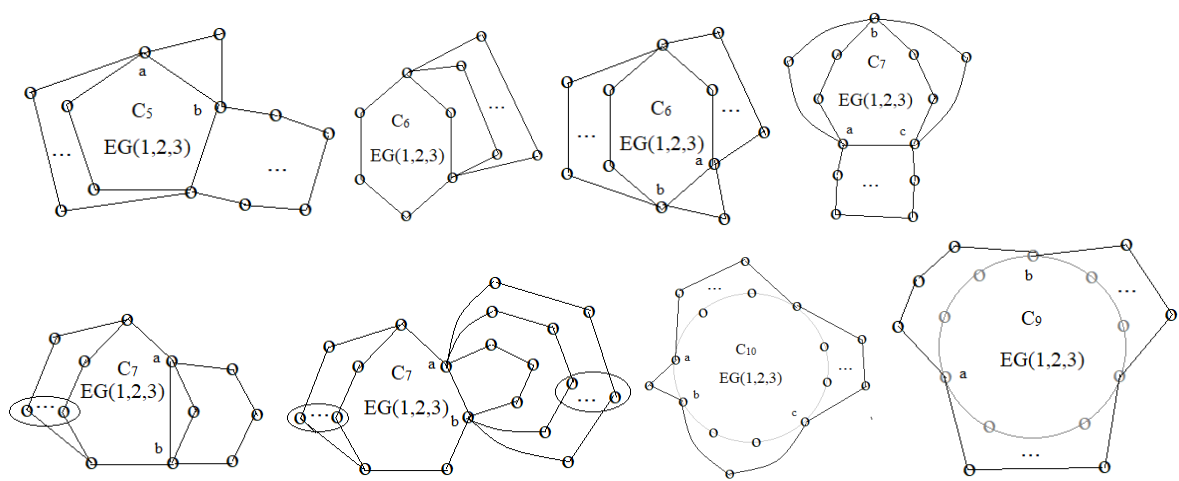

Fig.13 Infinite classes of graphs from $\varepsilon_{123}$.

Cycle intersection (Combined cycle) rules for Euler graphs from $\varepsilon_{123}$ are given in the Table-11. Graph in Fig.13a has a handle H(a,b) of length 2 and cannot be repeated, else cycle of type ≡0(mod 4) exists. Figs.13c,d,e,f,g,h also have non-repeatable handles.

**Appendix**

**Regular Graphs**

Markus Meringer, (2009) has considered enumeration of regular graphs for given number of vertices and degrees. Two tables of interest are:
- Connected regular graphs: The table contains numbers of connected regular graphs with 4-26 vertices and degree 3-7. The table indicates that the number of connected regular graphs of degree 2t is larger than that of the number of connected regular graphs of degree 2t-1. Also note that the number of connected regular graphs of odd degree of odd order is zero. It is well known that regular graphs of odd degree are of even order.
- Connected bipartite regular graphs: The table contains numbers of connected bipartite regular graphs with 6-32 even number of vertices and degree 3-5.

**Other Regularities**

Regularities with additional constraints were of interest in the literature, for instance, strongly regularity and distance regularity See Bose (1962), Cameron (2001), Paauwe (2007), Weisstein (2014). While 3-Cube is a distance regular graph, regular graphs defined from Partial Geometries, Partially Balanced Block designs and Mutually Orthogonal Latin Squares are strongly regular graphs. Existence and Construction of these graphs in general is still a challenge. Block designs are the backbone structures in the Experimental Designs with many applications: Agricultural Experiments, Error Correcting Codes, controlled sampling, randomized response, validation and valuation studies, intercropping experiments, brand cross-effect designs. Regularity attribute in natural or manmade systems greatly simplifies its makeup and leads to enhanced understanding of the underlying phenomenon.

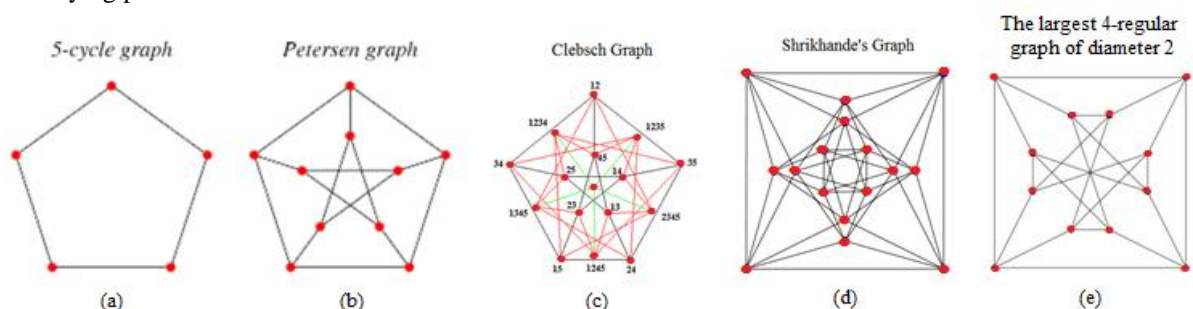

Fig.14a,b,c,d Strongly Regular Graphs, Fig.14e. Largest 4-regular graph of diameter 2.

Fractal world demonstrated that complex patterns of regular structures are feasible. Fractals having exact self similarity are called *regular fractals* (P. Addison (1997)). Regular Fractals are structures comprising of exact copies of themselves at all magnifications. For example, Cantor set.

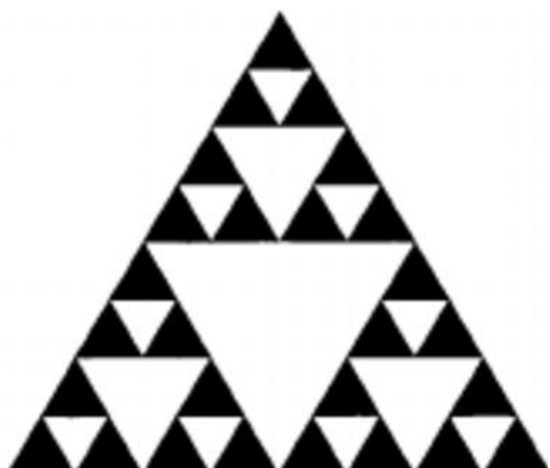

Fig.15 The Sierpinksi Gasket

*Regular fractal*. In regular fractals, the same geometric shape is repeated at multiple scales. The Sierpinksi Gasket is a typical regular fractal.

For six Regularities in the Solar System see Sherwood Harrington (2016).

**Evenness**

Evenness enters naturally at parties with dancing couples or handshaking tradition. Graph Theory presents them as Theorems. Euler (1736) proved that in a finite, undirected graph, *sum of node degrees is always even and equals to twice the number of edges* (Handshaking Lemma). A consequence is: *in any graph number of nodes of odd degree is even*. Two other instances evenness dominates in Euler graphs: 1. Number of edge disjoint paths between any two nodes of Euler graph is even. 2. A connected graph G is Euler graph if and only if all nodes of G are of even degree.

Regular Euler graphs with only one type or only three types of cycles are nonexistent or rare. Regularity in the class of Euler graphs with only two types of cycles, viz., (0,1), (0,3), (1,2), (1,3), (2,3) is nonexistent. It is conjectured that with only three types of cycles viz., (0,1,2), (0,1,3), (0,2,3), (1,2,3) regularity is absent except for $K_5$. However, Euler graphs are rich in regularity in the case (0,2), that is, when the graphs are bipartite Euler. Note that in this case all cycles are also even. Regular graphs are numerous once all four types of cycles are allowed.

*Cosmos within Chaos* is well deliberated and understood in the past. It is apparent that irregularity is mostly pervaded in the universe. However, regularity is an attribute which has been alluring, useful and is of great interest to humanity. Regularity is built-in or is a result of complex processes within the Universal systems. As the cosmic clouds condense, beautiful regular structures start forming. Star Systems evolve with regularity in the vast Universe. Our Solar system shows regularity in many respects. The regularities may be graded based on complexity. For example, complex regularities work within Earth systems.

The significance of evenness and regularity need further investigation.

## Acknowledgements

The author, an Oil & Gas Professional with works in Graph Theory, expresses deep felt gratitude to
The *Sonangol Pesquisa e Produção*, Luanda, Angola with special thanks to peers
*Joao Noguiera, Alexander Rosa, Solomon W Golumb,*
*T.V. Seshagiri Rao, J.L. Narasimham, G.A. Patwardhan and L. Radhakrishana*
and
The family *Usha, India; Chandana & Madhu, London, UK;*
*Neha & Avinash, Anika, Kiaan, Ashburn, Virginia,USA; Chetan & Vimmi,Tanush, SunnyVale, California, USA.*